%% file: hamBrieOrder.tex
\documentclass[11pt]{article}

\PassOptionsToPackage{rgb,x11names,svgnames,dvipsnames}{xcolor}
\usepackage[usenames,dvipsnames]{xcolor}
\definecolor{burnum}{rgb}{0.54, 0.2, 0.14}
\usepackage{tikz}
\usetikzlibrary{arrows}
\usetikzlibrary{shapes}
\usetikzlibrary{decorations.pathreplacing,decorations.markings}
\tikzset{->-/.style={decoration={
  markings,
  mark=at position .5 with {\arrow{>}}},postaction={decorate}}}
\usepackage{enumerate}
\usepackage{caption}
\usepackage{subcaption}

\input pre.tex

\begin{document}

\Large

\centerline{\bf{The odd-even invariant}}
\centerline{\bf{ and}}
\centerline{\bf{ Hamiltonian circuits in tope graphs}}
\smallskip
\centerline{by}
\centerline{Yvonne Kemper\footnote{Supported in part by Austrian Science Foundation FWF (START grant Y463).  Part of this research was performed while this author held a National Institute of Standards and Technology National Research Council Postdoctoral Associateship at the Information Technology Lab.} and Jim Lawrence}
\bigskip

\abstra{
In this paper we consider the question of the existence of
Hamiltonian circuits in the tope graphs of central arrangements of hyperplanes.
Some of the results describe connections
between the existence of Hamiltonian circuits in the arrangement
and the odd-even invariant
of the arrangement.  In conjunction with this, we present some results
concerning bounds on the odd-even invariant.
The results given here can be formulated more generally for oriented matroids and are still
valid in that setting.}

\sectb\sa. Introduction.  An arrangement $\A$ of
hyperplanes $H_i$, $i=1,\ldots,n$, in $\real^d$ determines a partition of $\real^d$ into
relatively open convex polyhedra.  Following \cite{ombook}, we call
the $d$-dimensional cells of this partition the {\it topes} of $\A$, and denote this set by $\T(\A)$.
In other words, the topes of the arrangement are the connected components of
the complement of the union $\bigcup_iH_i$.
The {\it tope graph} of $\A$ is the graph whose vertex set is $\T(\A)$,
with two topes (vertices) being adjacent if the intersection of their
closures is a facet of each.
The tope graph is always a connected bipartite graph.  There is thus a coloring of the topes by two colors -- say, burnt umber and chartreuse -- such that no two adjacent topes receive the same color, and this coloring
is unique up to reversing the colors of all the topes.  The question of which pairs $(b,c)$ are possible for given $n$ and $d$, where $b$ and $c$ are, respectively, the number of burnt umber and chartreuse topes, has been studied since at least the 1970s.  See, for example, \cite{grunb, harb, palasti, purdywetz, simmons, simmwetz}.

It is easy to see that if the tope graph possesses a Hamiltonian circuit,
then $\lvert b-c \rvert$ (the odd-even invariant) is 0.  This paper was largely
 motivated by a question of Vic Reiner \cite{reiner}, who asked if the
converse holds.
This question is studied here for central arrangements, that is, arrangements in
which all the hyperplanes contain the origin.
However the paper also includes results on bounds for the odd-even invariant concerning
the broader class of (not necessarily central) arrangements.
  In Section \sb\ we give the relevant definitions and background as well as several examples.  In Section \sd\ we discuss the specific case of alternating arrangements and Hamiltonian circuits, and in Section \se, we
establish the existence of arrangements having large odd-even invariants.
In Section \sf,
Theorem \th\ establishes the existence of arbitrarily large central arrangements in $\real^d$
having odd-even invariant 0 and having no Hamiltonian circuit.  Theorem \tj\ yields the fact
 that when $d$ is odd, the limit superior
as $n$ goes to infinity of $\frac b c$ for simple arrangements (arrangements in affine
general position) of $n$ hyperplanes
in $\real^d$ is 1.  This settles a question mentioned by Gr\"unbaum in \cite{grunb} and studied extensively when
$d=3$
by Purdy and Wetzel in \cite{purdywetz}.
Theorem \tk\ gives a criterion to verify that certain tope graphs do not possess perfect matchings.

Although the paper technically answers Reiner's question, it represents only a partial answer
to the issue raised of the relationship between the odd-even invariant and the existence
of Hamiltonian circuits in the tope graph.
The paper concludes in Section \sg\ with several
questions of interest.

The results are
presented in the setting of hyperplane arrangements in $\real^d$;
however, the results obtain also for oriented matroids and the arguments
can easily be adapted to the more general setting by appeal to the
topological representation theorem for oriented matroids.
See \cite{ombook} for a statement of this theorem.

\sectb\sb. Background and examples.
Let $\A$ denote an {\it arrangement of hyperplanes} in $\real^d$; that is, $\A$ is a finite, indexed collection of
 ($d-1$)-dimensional affine
subspaces of $\real^d$.  If $n$ denotes the number of hyperplanes, where $n\ge0$,
we may write $\A = \{H_1, H_2, \ldots, H_n\}$.

If the hyperplanes $H_i$ are linear subspaces the arrangement is said to be {\it central}.  Further, an arrangement of hyperplanes in $\real^d$ is {\it simple} if for $0\le k\le d$
the intersection of each set of $k$ of the hyperplanes of the arrangement has dimension $d-k$,
and the intersection of any $d+1$ or more of the hyperplanes is empty.  Equivalently, the
hyperplanes of the arrangement are in affine general position.
A central arrangement in $\real^d$ having more than $d$ hyperplanes cannot be simple.
A central arrangement of hyperplanes in $\real^d$ is {\it centrally simple} provided
that for $0\leq k \le d$ the intersection of each $k$ of the hyperplanes has
dimension $d-k$.  Equivalently, the hyperplanes of the arrangement are in linear general position.

The {\it rank} of the central arrangement
$\A$ is the difference $d-w$, where $w$ is the dimension of the
linear space that is the intersection of the $n$ hyperplanes.
We will usually assume that $w = 0$ (that is, the intersection contains only the origin), so that the rank of $\A$ is $d$.  If $V$ is a linear subspace
that is complementary to the linear subspace
$\bigcap_i H_i$, then the combinatorial structure of the
arrangement is adequately reflected in the arrangement $\{H_i^\prime = H_i\cap V\}$ of
hyperplanes in $V$, and for this arrangement, the dimension of $V$ is
the rank.

We will assume that positive and negative sides of each hyperplane in $\A$ are specified.
The open halfspace bounded by $H$ and lying on its positive side will be denoted $H^+$, the other $H^-$.
A tope $T$ may be specified by a function $s_T:\A\rightarrow\{-1,1\}$, where
$s_T(H) = 1$ if $T\subseteq H^+$ and $s_T(H) = -1$ if $T\subseteq H^-$.
Then $T = \bigl(\bigcap_{H:s_T(H) = 1}H^+\bigr)\cap(\bigcap_{H:s_T(H) = -1}H^-\bigr)$.

The collection of all possible functions $s:\A\rightarrow\{-1,1\}$ forms the vertex set
of the {\it cube graph} of order $n$, with two such functions (vertices)
 being adjacent if their values
 differ on exactly one hyperplane $H\in\A$.
In the situation of most interest in this paper, the (indexed) hyperplanes
of $\A$ are distinct: if $i\ne j$ then $H_i \ne H_j$.  In this case, the tope graph is a subgraph of the cube graph.  The cube graph is bipartite, and we may $2$-color its vertices with the colors 
chartreuse and burnt umber by assigning the color burnt umber to
any vertex $s$ (representing the map $s:\A\rightarrow \{-1,1\}$) for which the set
$\{H\in \A : s(H) = 1\}$ has even cardinality, and otherwise assigning
the color chartreuse to $s$.  Then, if $s_1$ and $s_2$ are adjacent, they have
different colors.
When the hyperplanes of $\A$ are distinct, the tope graph is
 a subgraph of the cube graph, and the 2-coloring of the cube graph yields a
2-coloring of the tope graph.

The {\it odd-even invariant} (as defined in \cite{refoea})
is 
\[\oe(\A) = \lvert \sum_{T\in\T(\A)}(-1)^{\sigma(T)}\rvert,\]
where $\sigma(T)$ denotes the number of hyperplanes $H$ such that
$T\subseteq H^+$.  Therefore, $\oe(\A)$ is the absolute value of the difference
between the number of burnt umber topes and the number of chartreuse topes.
We can further define the {\it signed odd-even invariant} as
\[\soe(\A) = \sum_{T\in\T(\A)}(-1)^{\sigma(T)}.\]
This sum yields the number of burnt umber topes minus the number of chartreuse topes.

Any central arrangement of $n$ distinct hyperplanes, with $n$ odd, has
odd-even invariant zero.  To see this, note that $T = \bigl(\bigcap_{H:s_T(H) = 1}H^+\bigr)\cap(\bigcap_{H:s_T(H) = -1}H^-\bigr)$ is a
 tope if and only if $T^*$ is a tope, where $s_{T^*}(H)=-s_{T}(H)$ for all $H\in\A$.
Then $\sigma(T)=n-\sigma(T^*)$.  When we have an odd number of hyperplanes,
 $(-1)^{\sigma(T)}$ and $(-1)^{\sigma(T^*)}$ cancel out in the sum.

\thrm\ \ta. For any arrangement $\A$ whose tope graph admits a Hamiltonian circuit, $\oe(\A) = 0$.

\proof
Given a bipartite graph with a $2$-coloring (chartreuse and burnt umber) of the vertices, the colors along any Hamiltonian circuit must alternate, as no similarly colored vertices are adjacent.  Since every vertex is included, the number of chartreuse and burnt umber vertices must be equal. ~\endproof

The Hamiltonian circuits of certain special tope graphs have been well studied.  We give a few examples here.\\
\noindent\textbf{Cube graphs.}  The cube graphs themselves are the tope graphs of the arrangements consisting
of the $n$ coordinate hyperplanes in $\real^n$.  The odd-even invariant of such arrangements is zero, and there exist Hamiltonian circuits of their tope graphs.  These circuits correspond to Gray codes (see \cite{stantonwhite}).\\
\noindent\textbf{Coxeter arrangements.} Another family of examples
comes from the Coxeter arrangements of type $A$.  The $\binom{n+1}{2}$ hyperplanes
of this arrangement are the sets $H_{i,j} = \{ (x_0, x_1, \ldots, x_{n}) \in V:
x_i = x_j\}$, where $V$ is the $n$-dimensional subspace of $\real^{n+1}$ satisfying
$\sum x_i = 0$.  A point $(x_0,x_1, \ldots, x_n) \in V$ that is not on
any of the hyperplanes has no two coordinates equal, so the coordinates are
ordered, $x_{\pi(0)}< x_{\pi(1)}< \ldots< x_{\pi(n)}$, for some permutation $\pi$.
Thus, the topes of the arrangement correspond to the permutations of $\{0,1,\ldots,n\}$.
Two topes are adjacent if and only if the corresponding permutations $\pi_1$ and
$\pi_2$ differ by a transposition $\tau$ that switches $i$ and $i-1$ for some $i$
with $1\leq i\leq n$.  That is, for $k=0,\ldots,n$, we have $\pi_1(\tau(k))=\pi_2(k)$, where
\[\tau(k) = \begin{cases} i& \text{if } k = i-1,\\
i-1& \text{if } k = i,\text{ and}\\
k&\text{otherwise.}
\end{cases}\]
These tope graphs are also well-known as the graphs of the permutahedra.  For these arrangements, the odd-even invariant is zero: the permutations of $\{0,\ldots, n\}$
are divided into the two color classes by assigning the even permutations one
of the colors and the odd permutations the other.
Hamiltonian circuits for general $n$ are given, for example, in \cite{stantonwhite}.

\begin{figure}
\centering
\begin{tikzpicture}[scale=.8]
\draw[thick,gray, <->] (0,3) -- (0,-3);
\draw[ thick,gray, <->] (-3,2) -- (3,-2);
\draw[ thick,gray, <->] (-3,-2) -- (3,2);
\draw[very thick] (2,0) -- (1,1.5) -- (-1,1.5) -- (-2,0) -- (-1,-1.5) -- (1,-1.5) -- (2,0);
\filldraw[fill=burnum] (2,0) circle (0.1);
\filldraw[fill=Chartreuse] (-2,0) circle (0.1);
\filldraw[fill=Chartreuse] (1,1.5) circle (0.1);
\filldraw[fill=burnum] (-1,1.5) circle (0.1);
\filldraw[fill=burnum] (-1,-1.5) circle (0.1);
\filldraw[fill=Chartreuse] (1,-1.5) circle (0.1);
\node at (3,-2.25) {{\footnotesize $x_0=x_1$}};
\node at (3,2.2) {{\footnotesize $x_1=x_2$}};
\node at (0,3.2) {{\footnotesize $x_0=x_2$}};
\node at (0,-3.3) {{\footnotesize The tope graph of the}};
\node at (0,-3.7) {{\footnotesize Coxeter arrangement $A_2$.}};

\begin{scope}[xshift=6cm, yshift=-1cm]
\draw[gray, thick,dotted] (2.5,2.5) -- (2.5,.5) -- (.5,.5);
\draw[gray, thick,dotted] (0,0) -- (0,2) -- (2,2);
\draw[gray, thick,dotted] (0,2) -- (.5,2.5);
\draw[gray, thick,dotted] (2,2) -- (2.5,2.5);
\draw[gray, thick,dotted] (2,0) -- (2.5,.5);
\draw[thick] (0,0) -- (2,0) -- (2,2) -- (2.5,2.5) -- (.5,2.5);
\draw[thick,color=black!60!white] (.5,2.5) -- (.5,.5) -- (0,0);
\filldraw[fill=burnum] (0,0) circle (0.1);
\filldraw[fill=Chartreuse] (2,0) circle (0.1);
\filldraw[fill=Chartreuse] (.5,.5) circle (0.1);
\filldraw[fill=burnum] (.5,2.5) circle (0.1);
\filldraw[fill=burnum] (2,2) circle (0.1);
\filldraw[fill=Chartreuse] (2.5,2.5) circle (0.1);
\node at (1.25,-.5) {{\footnotesize The tope graph as a subgraph}};
\node at (1.25,-.9) {{\footnotesize of the cube graph $Q_3$.}};
\end{scope}
\end{tikzpicture}
\caption{A Hamiltonian circuit on the Coxeter arrangement $A_2$.}
\label{Cox2}
\end{figure}
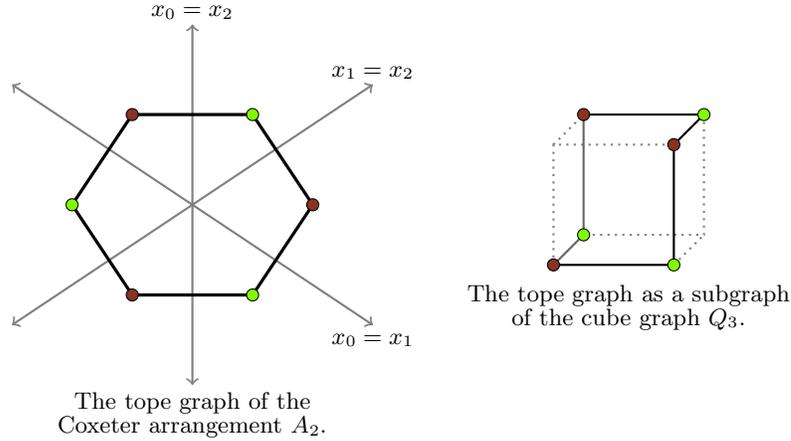

As a specific example, in Figure \ref{Cox2} we take $n=2$, with the hyperplanes $\{H_{01},H_{12},H_{02}\}=\{x_0=x_1,x_1=x_2,x_0=x_2\}$, shown projected on the subspace $x_0+x_1+x_2=0$.  (In this case, the Hamiltonian circuit is the entire tope graph.)

The tope graphs of the other finite Coxeter arrangements also admit Hamiltonian circuits.  See \cite{conslowil} for further details.  The next class of arrangements warrants its own section.\\
\sectb\sd. Alternating arrangements.
The following construction yields an
{\it  alternating
arrangement} of $n$ hyperplanes of rank $d$, and each combinatorial type of alternating
arrangement arises in this way.
  For $\alpha\in\real$,
let $H(\alpha) = \{(x_0,x_1,\ldots,x_{d-1}) : \sum_{k=0}^{d-1} x_k\alpha^k = 0\}$.
Given real numbers $\alpha_1<\alpha_2<\cdots<\alpha_n$,
let $H_i = H(\alpha_i)$ for $1\le i\le n$.
The arrangement $\A(n,d) = \{H_i : 1\le i \le n\}$ is an alternating
arrangement.

The topes of the arrangement are given in the following way.  Each vector $x=(x_0,x_1,\ldots,x_{d-1})\in\real^d$
gives rise to a real polynomial $p$ of degree $d-1$  whose
coefficients are the entries of $x$:
$$p(\tau) = x_0+x_1\tau+x_2\tau^2+\ldots+x_{d-1}\tau^{d-1}.$$
We may identify $\real^d$ with the vector space of such polynomials.
Two polynomials $p$, $q$ lie in the same tope of $\A$ if
$p$ and $q$ have the same nonzero sign when evaluated at each of the $\alpha_i$'s.

Alternatively, a function $s:\A\rightarrow\{-1,1\}$ corresponds to some tope $T$
if and only if the sequence $(s(H_1), s(H_2), \ldots, s(H_n))$
has at most $d-1$ sign changes, for if this is the case then there
is a polynomial of degree at most $d-1$ having the required sign
pattern.  The combinatorial type therefore does not depend upon the particular
choice of $\alpha_i$'s, but only upon their order.  The odd-even invariants of the alternating arrangements were determined
in \cite{refoea}; we restate the relevant theorems here.  

 \thrm\ \tb.  The odd-even invariant of the alternating arrangement
$\A(n,d)$, with $n$ even and $d$ odd, is
$2\genfrac{(}{)}{0pt}{0}{\frac{n}{2} -1}{\frac{d-1}{2}}$.

For the remaining cases, that is, when $n$ is odd and $d$ is even, recall that for any central
arrangement the odd-even invariant
is always zero when $n$ is odd.  The same is true when
$d$ is even, provided the arrangement is centrally simple
(as are the alternating arrangements). 

\thrm\ \tc.  If $\A$ is a centrally simple arrangement of $n$ hyperplanes in $\real^d$, then $\oe(\A) = 0$ unless $n$ is even and $d$ is odd.

Considering the notion of a ``mutation'' lends some credence to Theorem \tc, and this notion could be used to provide a proof that is valid in the current setting.  (The proof in \cite{refoea} does not use
mutations, but is valid for uniform oriented matroids in general.)

\begin{figure}
\centering
\begin{subfigure}[b]{0.4\textwidth}
\begin{center}
\begin{tikzpicture}
\filldraw[color=Chartreuse,draw=none] (.45,-.45) -- (0,0) -- (.45,.45) -- (.45,-.45);
\draw[thick] (1,1) -- (-1,-1);
\draw[thick] (-1,1) -- (1,-1);
\draw[thick] (.45,1.35) -- (.45,-1.25);
\draw[thick,dotted,color=red,->] (.425,.75) -- (-.75,.5);
\draw[thick,dotted,color=red,->] (.425,-.75) -- (-.75,-.5);
\node at (.315,1.3) {\tiny{$+$}};
\node at (.585,1.3) {\tiny{$-$}};
\node at (.85,1.05) {\tiny{$+$}};
\node at (1.05,.85) {\tiny{$-$}};
\node at (.85,-1.05) {\tiny{$-$}};
\node at (1.05,-.85) {\tiny{$+$}};
\end{tikzpicture}
\end{center}
\caption{}
\label{MutA}
\end{subfigure} \begin{subfigure}[b]{0.4\textwidth}
\begin{center}
\begin{tikzpicture}
\filldraw[color=burnum,draw=none] (-.45,-.45) -- (0,0) -- (-.45,.45) -- (-.45,-.45);
\draw[thick] (1,1) -- (-1,-1);
\draw[thick] (-1,1) -- (1,-1);
\draw[thick] (-.45,1.35) -- (-.45,-1.25);
\node at (-.315,1.3) {\tiny{$-$}};
\node at (-.585,1.3) {\tiny{$+$}};
\node at (.85,1.05) {\tiny{$+$}};
\node at (1.05,.85) {\tiny{$-$}};
\node at (.85,-1.05) {\tiny{$-$}};
\node at (1.05,-.85) {\tiny{$+$}};
\end{tikzpicture}
\end{center}
\caption{}
\label{MutB}
\end{subfigure}
\caption{Mutating a triangular tope in a planar arrangement.}
\end{figure}
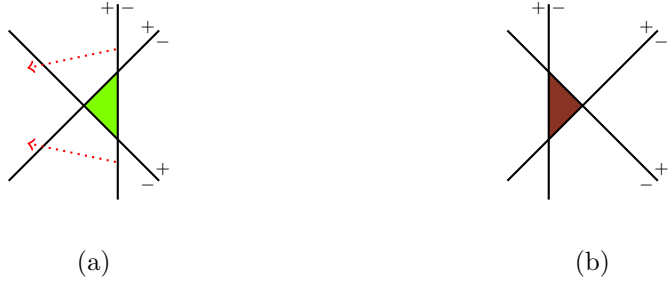

Given an arrangement in $\real^2$ that contains the three lines pictured in Figure \ref{MutA}, a {\it mutation} of the arrangement involving the triangular tope bounded by the three lines leads to an arrangement partly pictured in Figure \ref{MutB} (we only mutate simplicial topes).  Note that there may be more lines in the arrangement, but we may assume them to be unaffected by the mutation.

When positive sides of the lines are chosen,
the parity, even or odd, of the number of negative sides that contain
the tope changes, since the new tope lies on different sides of the
three lines from the old tope.  Therefore, we increase the number of one color tope by one, and decrease the other by one (in this case, we lose a chartreuse tope, and gain a burnt umber tope, but all else remains the same). 

In the analogous situation involving central arrangements,
two opposite topes, each an open simplicial cone, are replaced by
two other simplicial topes.  When $n$ is odd, a tope and its opposite have different colors to begin with, thus the odd-even invariant is unchanged.  When $n$ is even and $d$ is odd, the odd-even invariant changes by $\pm4$, as we change the signs of an odd number of hyperplanes for each tope (that is, the simplicial cone is bounded by an odd number of hyperplanes, and we change the sign of all the bounding hyperplanes in $s_T$).  However, when $d$ is even, the odd-even invariant does not change, as we switch the signs of an even number of hyperplanes for each antipodal tope.

For certain pairs $(n,d)$, we are able to describe Hamiltonian circuits in the tope graphs of the alternating arrangements $\A(n,d)$.

\thrm\ \td.
The tope graph of the alternating arrangement $A(n,3)$ 
with $n$ odd is Hamiltonian.

\proof Consider the alternating arrangement as the set of sequences of $\{+,-\}$'s of length $n$ and fewer than $3$ changes of sign.  Given a sequence $s=s_1s_2\cdots s_n\in\{+,-\}^n$, we
 may represent the sequence as the set $S=\{i :s_i=-\}$.  When
 we allow at most two sign changes, this means that the $i\in
 S$ are consecutive modulo $n$.
We denote the possible sets $S$ by $S_0$, $S_n$, and, for $1\leq j\leq n-1$
and $1\le k\le n$, $S_{j,k}$, described as follows.
We put $S_0=\emptyset$, and $S_n=\{1,2,\ldots,n\}$. 
 For $1\leq j\leq n-1$ and $1\leq k\leq n$, let $S_{j,k}$ be the set of $j$ consecutive integers (modulo $n$)
starting at
 position $k$,  for instance,
 $S_{1,3} = \{3\}$ and $S_{3,n} = \{n,1,2\}$.
These $2 + n(n-1)$ sets correspond to the topes, which are the vertices of the tope graph.
We utilize the sets to describe the tope graph.  The graph has two vertices of degree $n$,
$n(n-3)$ vertices of degree 4, and $2n$ vertices of degree 3.
The vertex $S_0$ is adjacent to each vertex $S_{1,k}$, and $S_n$
 is adjacent to $S_{n-1,k}$, both for all $k$, $1\leq k\leq n$;
$S_0$ and $S_n$ are the vertices of degree $n$.  When $2\le j\le n-2$ and
$1\le k\le n$,  $S_{j,k}$
 is adjacent to $S_{j-1,k}$, $S_{j-1,k+1}$, $S_{j+1,k}$, and
 $S_{j+1,k-1}$; these are the vertices of degree 4. 
(We consider $k$ modulo $n$; e.g. $S_{2,1}$ is adjacent to $S_{1,1}$, $S_{1,2}$, $S_{3,1}$ and $S_{3,0} = S_{3,n}$.)
Also, $S_{1,k}$ is adjacent to $S_{2,k}$ and $S_{2,k-1}$, in addition to $S_0$; and $S_{n-1,k}$ is adjacent to
$S_{n-1,k}$ and $S_{n-1,k+1}$, in addition to $S_n$.
 
The next portion of the proof involves a great many indices, and the following is (hopefully) the clearest presentation.  We populate an $n\times (n-1)$ array with the sets $S_{j,k}$ as follows.  Express $j$ ($1\le j \le n-1$)
as $j=4t+r$, where $t$ and $r$ are integers, and, as $n$ is odd, $r\in\{1,3\}$.  Then, $S_{j,k}$ is placed in the $j$-th column and the $(k+t)$-th row.  It is harmless to assume (this will be justified later) that $n$ is of the form $4s+3$, in which case the array looks something like:
 
 {{\small
 \[\begin{array}{|c|c|c|c|c|c|c|c|c|c|}
 \hline
 S_{1,1} & S_{2,1} & S_{3,1} & S_{4,n} & S_{5,n} & \cdots & S_{n-4,2-s} & S_{n-3,1-s} & S_{n-2,1-s} & S_{n-1,1-s} \\
  \hline
S_{1,2} & S_{2,2} & S_{3,2} & S_{4,1} & S_{5,1} & \cdots & S_{n-4,3-s} & S_{n-3,2-s} & S_{n-2,2-s} & S_{n-1,2-s} \\
 \hline
S_{1,3} & S_{2,3} & S_{3,3} & S_{4,2} & S_{5,2} & \cdots & S_{n-4,4-s} & S_{n-3,3-s} & S_{n-2,3-s} & S_{n-1,3-s} \\
 \hline
\vdots & \vdots & \vdots & \vdots & \vdots &  & \vdots & \vdots & \vdots & \vdots \\
 \hline
S_{1,n-1} & S_{2,n-1} & S_{3,n-1} & S_{4,n-2} & S_{5,n-2} & \cdots & S_{n-4,n-s} & S_{n-3,n-1-s} & S_{n-2,n-1-s} & S_{n-1,n-1-s} \\
 \hline
S_{1,n} & S_{2,n} & S_{3,n} & S_{4,n-1} & S_{5,n-1} & \cdots & S_{n-4,n+1-s} & S_{n-3,n-s} & S_{n-2,n-s} & S_{n-1,n-s} \\
 \hline
 \end{array}\]
 }}
 
\noindent The Hamiltonian circuit we describe begins at $S_0$, which is adjacent to every vertex in the first column of our array.  From there, we go to array element $(1,1)$, and then zig-zag between pairs of columns, leaving one element in each column unvisited; these will be used to travel back from $S_n$.  In particular, for columns $y$ with $y\equiv 1\bmod 4$ and rows $x$, $1\leq x\leq n-2$, we travel:
\[\cdots\rightarrow (x,y) \rightarrow (x,y+1) \rightarrow (x+1,y) \rightarrow (x+1,y+1) \rightarrow \cdots.\]
For $x=n-1$, we continue to the next pair of columns, zig-zag our way back up to the first row, and on to the next pair of columns:

\[\begin{array}{rlllllll}
\cdots\rightarrow & (n-1,y) & \rightarrow & (n-1,y+1) & \rightarrow & (n-1,y+2) & \rightarrow \\
 & (n-1,y+3) & \rightarrow & (n-2,y+2) & \rightarrow & (n-2,y+3) & \rightarrow\cdots \\
& & & \vdots & & & &\\
\cdots\rightarrow & (1,y+2) & \rightarrow & (1,y+3) & \rightarrow & (1,y+4) & \rightarrow \\
& (1,y+5) & \rightarrow & (2,y+4) & \rightarrow & (2,y+5) & \rightarrow\cdots\\
\end{array}\]
When we reach entry $(n-1,n-1)$ (in the case of $n\equiv 3\bmod 4$) or entry $(1,n-1)$ (if $n\equiv 1\bmod 4$), we continue to $S_n$.  (Every $S_{n-1,k}$ is adjacent to $S_n$, thus it does not matter if $n\equiv 1\bmod 4$ or if $n\equiv 3\bmod 4$ -- i.e. whether we zig-zag up or down the last pair of columns.)  Then, from $S_n$, we travel back to $S_0$:
\[S_n\rightarrow (n,n-1) \rightarrow (n,n-2) \rightarrow\cdots\rightarrow (n,2)\rightarrow (n,1)\rightarrow S_0,\]
completing the Hamiltonian circuit.  In diagram form (again, assuming $n=4s+3$), we have:
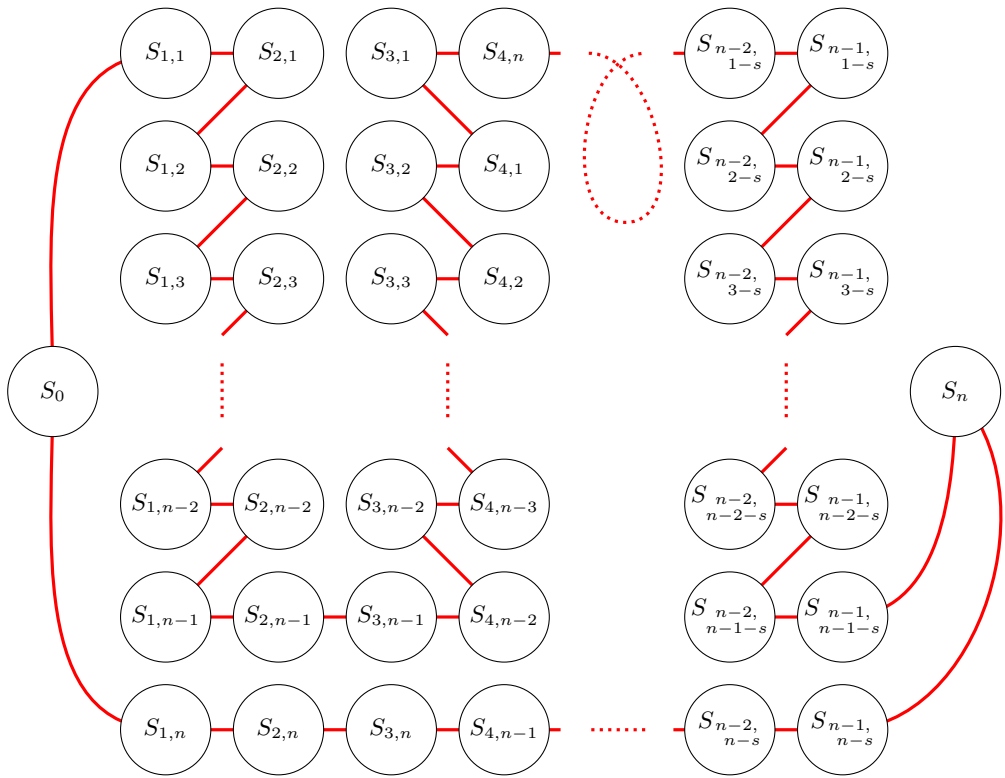
\begin{center}
\begin{tikzpicture}[scale=.75,font=\footnotesize]
\draw[red, very thick] (0,-8) to[out=90,in=180] (2,-2) -- (4,-2) -- (2,-4) -- (4,-4) -- (2,-6) -- (4,-6) -- (3,-7);
\draw[red, very thick, dotted] (3,-7.5) -- (3,-8.5);
\draw[red,very thick] (3,-9) -- (2,-10) -- (4,-10) -- (2,-12) -- (4,-12) -- (6,-12) -- (8,-12) -- (6,-10) -- (8,-10) -- (7,-9);
\draw[red, very thick, dotted] (7,-7.5) -- (7,-8.5);
\draw[red,very thick] (7,-7) -- (6,-6) -- (8,-6) -- (6,-4) -- (8,-4) -- (6,-2) -- (8,-2) -- (9,-2);
\draw[red,very thick,dotted] (9.5,-2) 
to[out=0,in=0] (10.15,-5) 
to[out=180,in=180] (10.5,-2);

\draw[red, very thick] (11,-2) -- (14,-2) -- (12,-4) -- (14,-4) -- (12,-6) -- (14,-6) -- (13,-7);
\draw[red, very thick, dotted] (13,-7.5) -- (13,-8.5);
\draw[red,very thick] (13,-9) -- (12,-10) -- (14,-10) -- (12,-12) -- (14,-12) to[out=0,in=270] (16,-8) to[out=-45,in=0] (14,-14) -- (12,-14) -- (11,-14);
\draw[red, very thick, dotted] (10.5,-14) -- (9.5,-14);
\draw[red, very thick] (9,-14) -- (2,-14) to[out=180,in=270] (0,-8);

\filldraw[fill=white] (0,-8) circle (.8);
\node at (0,-8) {$S_0$};
\filldraw[fill=white]  (16,-8) circle (.8);
\node at (16,-8) {$S_n$};
\foreach \x in {1,2,3,4,6,7}{
      \foreach \y in {-1,-2,-3,-5,-6,-7}{
        \filldraw[fill=white] (2*\x,2*\y) circle (.8);
}
}

\begin{scope}[scale=4/3]
\node at (1.5,-1.5) {$S_{1,1}$};
\node at (1.5,-3) {$S_{1,2}$};
\node at (1.5,-4.5) {$S_{1,3}$};
\node at (3,-1.5) {$S_{2,1}$};
\node at (3,-3) {$S_{2,2}$};
\node at (3,-4.5) {$S_{2,3}$};
\node at (4.5,-1.5) {$S_{3,1}$};
\node at (4.5,-3) {$S_{3,2}$};
\node at (4.5,-4.5) {$S_{3,3}$};
\node at (6,-1.5) {$S_{4,n}$};
\node at (6,-3) {$S_{4,1}$};
\node at (6,-4.5) {$S_{4,2}$};

\node at (1.5,-7.5) {$S_{1,n-2}$};
\node at (1.5,-9) {$S_{1,n-1}$};
\node at (1.5,-10.5) {$S_{1,n}$};
\node at (3,-7.5) {$S_{2,n-2}$};
\node at (3,-9) {$S_{2,n-1}$};
\node at (3,-10.5) {$S_{2,n}$};
\node at (4.5,-7.5) {$S_{3,n-2}$};
\node at (4.5,-9) {$S_{3,n-1}$};
\node at (4.5,-10.5) {$S_{3,n}$};
\node at (6,-7.5) {$S_{4,n-3}$};
\node at (6,-9) {$S_{4,n-2}$};
\node at (6,-10.5) {$S_{4,n-1}$};
\end{scope}

\begin{scope}[xshift=-2cm]
\begin{scope}[scale=4/3]
\node at (10.5,-1.5) {$S_{\substack{n-2,\\ \;\;\;1-s}}$};
\node at (10.5,-3) {$S_{\substack{n-2,\\ \;\;\;2-s}}$};
\node at (10.5,-4.5) {$S_{\substack{n-2,\\ \;\;\;3-s}}$};
\node at (12,-1.5) {$S_{\substack{n-1,\\ \;\;\;1-s}}$};
\node at (12,-3) {$S_{\substack{n-1,\\ \;\;\;2-s}}$};
\node at (12,-4.5) {$S_{\substack{n-1,\\ \;\;\;3-s}}$};

\node at (10.5,-7.5) {$S_{\substack{n-2,\\ n-2-s}}$};
\node at (10.5,-9) {$S_{\substack{n-2,\\ n-1-s}}$};
\node at (10.5,-10.5) {$S_{\substack{n-2,\\ \;\;n-s}}$};
\node at (12,-7.5) {$S_{\substack{n-1,\\ n-2-s}}$};
\node at (12,-9) {$S_{\substack{n-1,\\ n-1-s}}$};
\node at (12,-10.5) {$S_{\substack{n-1,\\ \;\;n-s}}$};
\end{scope}
\end{scope}
\end{tikzpicture}
\end{center}
~\endproof

As an example, we give the described Hamiltonian circuit for $\A(5,3)$:
{\small{
\[\begin{array}{l}
+++++,\ -++++,\ --+++,\ +-+++,\ +--++,\ ++-++,\\
++--+,\ +++-+,\ +++--,\ -++--,\ -+---,\ ++---,\\
+----, \ +---+,\ ----+,\ ---++,\ ---+-,\ -----,\\
--+--,\ --++-,\ -+++-,\ ++++-, \ +++++.\\
\end{array}\]}}

We may also view this circuit on the (spherical projection of the) arrangement itself; this is shown in Figure \ref{HamA53}.

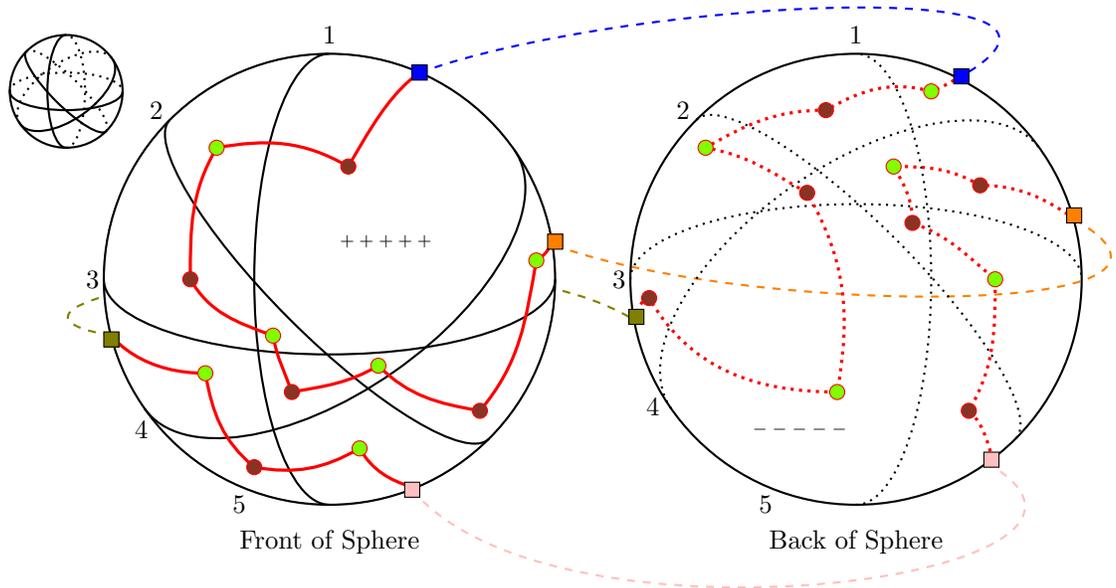
\begin{figure}
\centering
\begin{tikzpicture}
\begin{scope}[xshift=-3.5cm,yshift=2.5cm,scale=.25]
\draw[thick] (0,0) circle (3);
\draw[thick] (3,0) arc (0:-180:3 and 1);
\draw[thick] (0,3) arc (90:270:1 and 3);
\draw[thick,rotate=35] (3,0) arc (0:-180:3 and 1.5);
\draw[thick,rotate=-45] (3,0) arc (0:-180:3 and .75);
\begin{scope}[rotate=180]
\draw[thick,dotted] (3,0) arc (0:-180:3 and 1);
\draw[thick,dotted] (0,3) arc (90:270:1 and 3);
\draw[thick,dotted,rotate=35] (3,0) arc (0:-180:3 and 1.5);
\draw[thick,dotted,rotate=-45] (3,0) arc (0:-180:3 and .75);
\end{scope}
\end{scope}

\draw[pink,dashed,thick] (1.1,-2.8) to[out=-50,in=-30] (8.75,-2.45);
\draw[olive,dashed,thick] (-2.9,-.75) to[out=170,in=150] (3.98,-.5);
\draw[orange,dashed,thick] (3,.45) to[out=-20,in=-30] (9.9,.8);
\draw[blue, dashed, thick] (1.2,2.75) to[out=20,in=30] (8.4,2.7);
\filldraw[fill=white] (0,0) circle (3);
\draw[thick] (0,0) circle (3);
\draw[thick] (3,0) arc (0:-180:3 and 1);
\draw[thick] (0,3) arc (90:270:1 and 3);
\draw[thick,rotate=35] (3,0) arc (0:-180:3 and 1.5);
\draw[thick,rotate=-45] (3,0) arc (0:-180:3 and .75);
\begin{scope}
\clip (0,0) circle (3);
\draw[very thick, color=red] (2,3) to[out=180,in=60] (.25,1.5) 
to[out=150,in=10] (-1.5,1.75)
to[out=240,in=90] (-1.85,0) 
to[out=-60,in=160] (-.75,-.75) 
to[out=-80,in=110] (-.5,-1.5) 
to[out=10,in=210] (.65,-1.15) 
to[out=-45,in=170] (2,-1.75) 
to[out=50,in=260] (2.75,.25) 
to[out=30,in=270] (3,1);

\draw[very thick, color=red] (2,-3) to[out=160,in=-60] (.4,-2.25) 
to[out=210,in=-10] (-1,-2.5) 
to[out=140,in=-80] (-1.65,-1.25) 
to[out=180,in=290] (-3.2,-.25);
\end{scope}

\filldraw[fill=Chartreuse,draw=red] (.4,-2.25) circle (.1); 
\filldraw[fill=burnum,draw=red] (-1,-2.5) circle (.1); 
\filldraw[fill=Chartreuse,draw=red] (-1.65,-1.25) circle (.1); 
\filldraw[fill=burnum,draw=red] (-1.85,0) circle (.1); 
\filldraw[fill=Chartreuse,draw=red] (-1.5,1.75) circle (.1); 
\filldraw[fill=burnum,draw=red] (.25,1.5) circle (.1); 
\filldraw[fill=Chartreuse,draw=red] (-.75,-.75) circle (.1); 
\filldraw[fill=burnum,draw=red] (-.5,-1.5) circle (.1); 
\filldraw[fill=Chartreuse,draw=red] (.65,-1.15) circle (.1); 
\filldraw[fill=burnum,draw=red] (2,-1.75) circle (.1); 
\filldraw[fill=Chartreuse,draw=red] (2.75,.25) circle (.1); 

\begin{scope}[xshift=7cm,rotate=180]
\draw[thick] (0,0) circle (3);
\draw[thick,dotted] (3,0) arc (0:-180:3 and 1);
\draw[thick,dotted] (0,3) arc (90:270:1 and 3);
\draw[thick,dotted,rotate=35] (3,0) arc (0:-180:3 and 1.5);
\draw[thick,dotted,rotate=-45] (3,0) arc (0:-180:3 and .75);
\end{scope}

\begin{scope}[xshift=7cm]
\node at (-.75,-2) {{\tiny $-----$}};
\clip (0,0) circle (3);
\draw[red,very thick, dotted] (-3,-.5) to[out=10,in=180] (-2.75,-.25) 
to[out=-60,in=180] (-.25,-1.5) 
to[out=80,in=-60] (-.65,1.15) 
to[out=150,in=-10] (-2,1.75) 
to[out=30,in=180] (-.4,2.25) 
to[out=30,in=165] (1,2.5) 
to[out=20,in=160] (2,2.8);

\draw[red,very thick, dotted] (2.3,-3) to[out=200,in=-40] (1.5,-1.75) 
to[out=60,in=270] (1.85,0) 
to[out=140,in=-20] (.75,.75) 
to[out=90,in=-60] (.5,1.5) 
to[out=0,in=150] (1.65,1.25) 
to[out=0,in=145] (3.2,.6);
\end{scope}

\begin{scope}[xshift=7cm,rotate=180]
\filldraw[fill=burnum,draw=red] (.4,-2.25) circle (.1);
\filldraw[fill=Chartreuse,draw=red] (-1,-2.5) circle (.1);
\filldraw[fill=burnum,draw=red] (-1.65,-1.25) circle (.1);
\filldraw[fill=Chartreuse,draw=red] (-1.85,0) circle (.1); %
\filldraw[fill=burnum,draw=red] (-1.5,1.75) circle (.1); %
\filldraw[fill=Chartreuse,draw=red] (.25,1.5) circle (.1);%
\filldraw[fill=burnum,draw=red] (-.75,-.75) circle (.1); %
\filldraw[fill=Chartreuse,draw=red] (-.5,-1.5) circle (.1);%
\filldraw[fill=burnum,draw=red] (.65,-1.15) circle (.1); %
\filldraw[fill=Chartreuse,draw=red] (2,-1.75) circle (.1); %
\filldraw[fill=burnum,draw=red] (2.75,.25) circle (.1); %

\end{scope}

\begin{scope}[xshift=8.3cm,yshift=2.6cm]
\filldraw[fill=blue] (0,0) rectangle (.2,.2);
\end{scope}
\begin{scope}[xshift=1.1cm,yshift=2.65cm]
\filldraw[fill=blue] (0,0) rectangle (.2,.2);
\end{scope}

\begin{scope}[xshift=9.8cm,yshift=.75cm]
\filldraw[fill=orange] (0,0) rectangle (.2,.2);
\end{scope}
\begin{scope}[xshift=2.9cm,yshift=.4cm]
\filldraw[fill=orange] (0,0) rectangle (.2,.2);
\end{scope}

\begin{scope}[xshift=3.98cm,yshift=-.6cm]
\filldraw[fill=olive] (0,0) rectangle (.2,.2);
\end{scope}
\begin{scope}[xshift=-3cm,yshift=-.9cm]
\filldraw[fill=olive] (0,0) rectangle (.2,.2);
\end{scope}

\begin{scope}[xshift=8.7cm,yshift=-2.5cm]
\filldraw[fill=pink] (0,0) rectangle (.2,.2);
\end{scope}
\begin{scope}[xshift=1cm,yshift=-2.9cm]
\filldraw[fill=pink] (0,0) rectangle (.2,.2);
\end{scope}

\node at (0,3.25) {\small $1$};
\node at (-2.3,2.25) {\small $2$};
\node at (-3.15,0) {\small $3$};
\node at (-2.5,-2) {\small $4$};
\node at (-1.2,-3) {\small $5$};

\begin{scope}[xshift=7cm]
\node at (0,3.25) {\small $1$};
\node at (-2.3,2.25) {\small $2$};
\node at (-3.15,0) {\small $3$};
\node at (-2.7,-1.7) {\small $4$};
\node at (-1.2,-3) {\small $5$};
\end{scope}

\node at (0,-3.5) {\small{Front of Sphere}};
\node at (7,-3.5) {\small{Back of Sphere}};
\node at (.75,.5) {{\tiny $+++++$}};
\end{tikzpicture}
\vspace{-1.5cm}
\caption{A Hamiltonian circuit on $\A(5,3)$, illustrating Theorem \td.  The bounding edges of the front and back of the sphere are identified as indicated in the top left image.  Dashed lines connecting squares indicate that the path goes from the front of the sphere to the back (or vice-versa) at this point.}
\label{HamA53}
\end{figure}
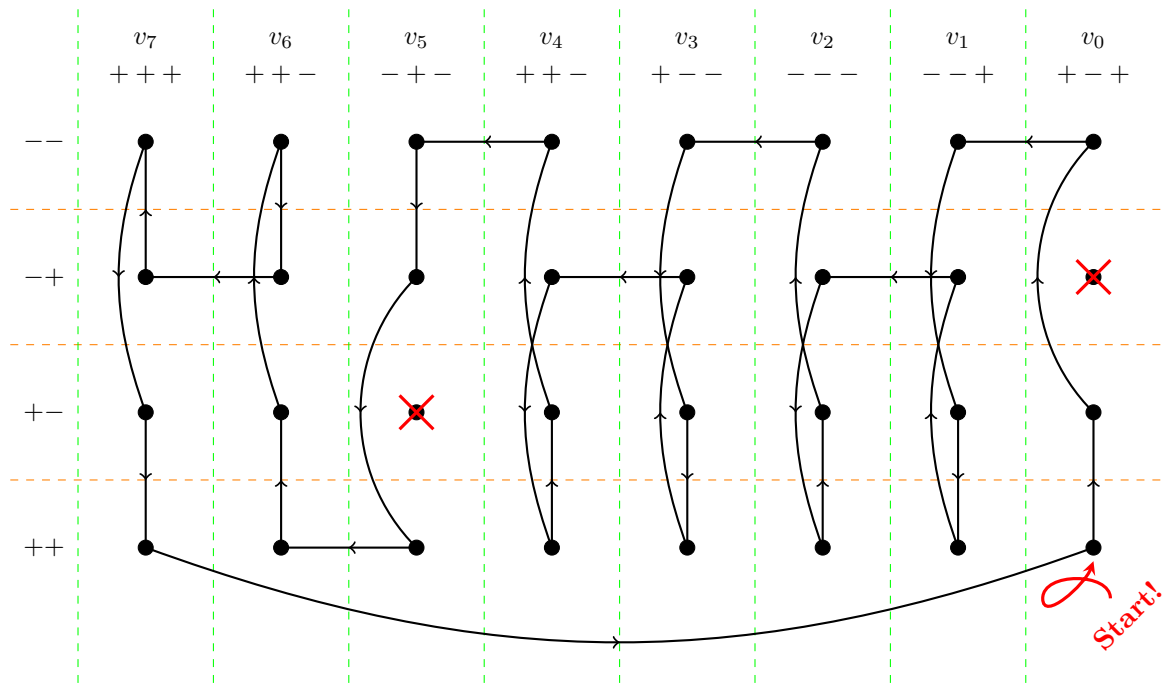
\begin{figure}
\begin{tikzpicture}[scale=.9,rotate=90]

\foreach \x in {1,2,3,4}{
      \foreach \y in {0,-1,-2,-3,-4,-5,-6,-7}{
        \node[draw,circle,inner sep=2pt,fill] at (2*\x,2*\y) {};
      }
}
\draw[very thick,color=red,->,>=stealth] (1.25,-14.25) to[out=0,in=0] (1.25,-13.25) to[out=180,in=160] (1.8,-14);
\node [rotate=45] at (1,-14.5) {{\color{red}\textbf{Start!}}};

\draw[dashed,color=orange] (3,2) -- (3,-15);
\draw[dashed,color=orange] (5,2) -- (5,-15);
\draw[dashed,color=orange] (7,2) -- (7,-15);

\draw[dashed,color=green] (0,1) -- (10,1);
\draw[dashed,color=green] (0,-1) -- (10,-1);
\draw[dashed,color=green] (0,-3) -- (10,-3);
\draw[dashed,color=green] (0,-5) -- (10,-5);
\draw[dashed,color=green] (0,-7) -- (10,-7);
\draw[dashed,color=green] (0,-9) -- (10,-9);
\draw[dashed,color=green] (0,-11) -- (10,-11);
\draw[dashed,color=green] (0,-13) -- (10,-13);

\node at (2,1.5) {\small{$++$}};
\node at (4,1.5) {\small{$+-$}};
\node at (6,1.5) {\small{$-+$}};
\node at (8,1.5) {\small{$--$}};

\node at (9,0) {\small{$+++$}};
\node at (9.5,0) {\small{$v_7$}};
\draw[->-,thick] (4,0) to (2,0);
\draw[->-,thick] (8,0) to[out=160,in=20] (4,0);
\draw[->-,thick] (6,0) to (8,0);
\node at (9,-2) {\small{$++-$}};
\node at (9.5,-2) {\small{$v_6$}};
\draw[->-,thick] (6,-2) to (6,0);
\draw[->-,thick] (8,-2) to (6,-2);
\draw[->-,thick] (4,-2) to[out=20,in=160] (8,-2);
\draw[->-,thick] (2,-2) to (4,-2);
\node at (9,-4) {\small{$-+-$}};
\node at (9.5,-4) {\small{$v_5$}};
\draw[->-,thick] (2,-4) to (2,-2);
\draw[very thick,red] (3.75,-3.75) -- (4.25,-4.25);
\draw[very thick,red] (3.75,-4.25) -- (4.25,-3.75);
\draw[->-,thick] (6,-4) to[out=135,in=45] (2,-4);
\draw[->-,thick] (8,-4) to (6,-4);
\node at (9,-6) {\small{$++-$}};
\node at (9.5,-6) {\small{$v_4$}};
\draw[->-,thick] (8,-6) to (8,-4);
\draw[->-,thick] (4,-6) to[out=20,in=160] (8,-6);
\draw[->-,thick] (2,-6) to (4,-6);
\draw[->-,thick] (6,-6) to[out=160,in=20] (2,-6);
\node at (9,-8) {\small{$+--$}};
\node at (9.5,-8) {\small{$v_3$}};
\draw[->-,thick] (6,-8) to (6,-6);
\draw[->-,thick] (2,-8) to[out=20,in=160] (6,-8);
\draw[->-,thick] (4,-8) to (2,-8);
\draw[->-,thick] (8,-8) to[out=160,in=20] (4,-8);
\node at (9,-10) {\small{$---$}};
\node at (9.5,-10) {\small{$v_2$}};
\draw[->-,thick] (8,-10) to (8,-8);
\draw[->-,thick] (4,-10) to[out=20,in=160] (8,-10);
\draw[->-,thick] (2,-10) to (4,-10);
\draw[->-,thick] (6,-10) to[out=160,in=20] (2,-10);
\node at (9,-12) {\small{$--+$}};
\node at (9.5,-12) {\small{$v_1$}};
\draw[->-,thick] (6,-12) to (6,-10);
\draw[->-,thick] (2,-12) to[out=20,in=160] (6,-12);
\draw[->-,thick] (4,-12) to (2,-12);
\draw[->-,thick] (8,-12) to[out=160,in=20] (4,-12);
\node at (9,-14) {\small{$+-+$}};
\node at (9.5,-14) {\small{$v_0$}};
\draw[->-,thick] (8,-14) to (8,-12);
\draw[->-,thick] (4,-14) to[out=45,in=135] (8,-14);
\draw[very thick,red] (5.75,-13.75) -- (6.25,-14.25);
\draw[very thick,red] (5.75,-14.25) -- (6.25,-13.75);
\draw[->-,thick] (2,-14) to (4,-14);
\draw[->-,thick]  (2,0) to[out=250,in=110] (2,-14);

\end{tikzpicture}
\caption{A path on $A(3,3)$ extended to $A(5,4)$, illustrating Theorem \tg.}
\label{A33toA54}
\end{figure}

If $n=d$, the tope graph of $A(n,n)$ is the cube
 graph of dimension $n$, which, as stated above, is also Hamiltonian. 
 We may use this fact to get a further class of $(n,d)$ pairs:

\thrm\ \tg. The tope graph of the alternating arrangement $\A(n,n-1)$, where $n$ is odd,
is Hamiltonian.

\proof
Set $N = 2^{n-2}-1$ and
suppose that $v_0, v_1, \ldots, v_{N}, v_0$ is a Hamiltonian circuit
in $\A(n-2,n-2)$, consecutive topes in this list being adjacent.
Suppose further that $v_0 = +-+-\cdots-+$ and $M$ is the integer, $0<M<N$, for
which $v_M = -+-+\cdots+-$.  Notice that $M$ must be odd.
Also note that neither $v_0-+$ nor $v_M+-$ is a tope in $\A(n,n-1)$; these are
the only two sequences of $n$ $+$'s and $-$'s that are not.

We describe a Hamiltonian circuit in $\A(n,n-1)$.
The path begins at $v_0++$ and proceeds to $v_0+-$ and then $v_0--$.
From $v_0--$ the circuit proceeds through $v_1--$, $v_1+-$, and $v_1++$
to $v_1-+$.
For  $k = 2, 3, \ldots, M-1$, when
$k$ is even,
the circuit proceeds from $v_{k-1}-+$ through $v_k-+$, $v_k++$, and $v_k+-$
to $v_k--$;
when $k$ is odd, the circuit proceeds from $v_{k-1}--$
through $v_k--$, $v_k+-$, and $v_k++$ to $v_k-+$.
After $v_{M-1}--$, the circuit passes through $v_M--$ and $v_M-+$ to $v_M++$.
From there it proceeds to $v_{M+1}++$, $v_{M+1}+-$, $v_{M+1}--$,
and $v_{M+1}-+$.
Then, for  $k = M+2, M+3,  \ldots, N-1$, when $k$ is odd, the circuit proceeds through
$v_k-+$, $v_k--$, $v_k+-$, and $v_k++$;
when $k$ is  even, it proceeds through $v_k++$,
$v_k+-$, $v_k--$, and $v_k-+$.  From $v_N++$ it returns to $v_0++$.
\endproof

We illustrate the case with $A(3,3)$ extending to $A(5,4)$ in Figure \ref{A33toA54}.  These initial results and observations, as well as some small examples, give rise to a very natural question: do the tope graphs of all alternating oriented matroids with odd-even invariant equal to zero have Hamiltonian circuits?

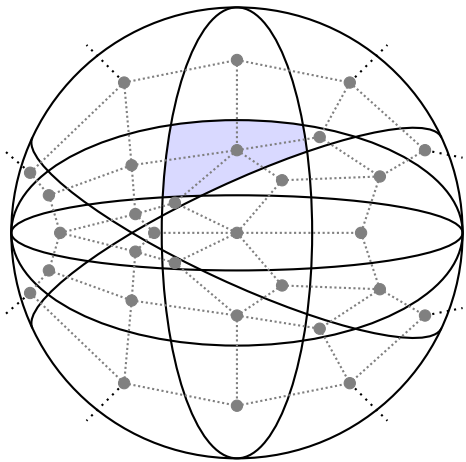
\begin{figure}[h]
\centering
\begin{tikzpicture}
\begin{scope} 
\path[clip] (3,0) arc (0:180:3 and 1.5) -- (-3,0) arc (180:0:3 and .5);
\path[clip] (0,0) ellipse (1 and 3);
\filldraw[fill=blue!15!white,draw=none,rotate=25] (3,0) arc (0:180:3 and .65) -- (-3,0) arc (180:0:3 and 3);
\end{scope}

\draw[thick] (0,0) circle (3);
\draw[thick] (0,0) ellipse (3 and .5);
\draw[thick] (0,0) ellipse (3 and 1.5);
\draw[thick] (0,0) ellipse (1 and 3);
\draw[thick,rotate=25] (3,0) arc (0:180:3 and .65);
\draw[thick,rotate=-25] (-3,0) arc (180:360:3 and .65);

\draw[thick,dotted] (-2.75,.8) -- (-3.1,1.1); 
\draw[thick,gray,densely dotted] (-2.75,.8) -- (-2.5,.5); 
\draw[thick,gray,densely dotted] (-2.75,.8) -- (-1.5,2); 
\draw[thick,gray,densely dotted] (-2.5,.5) -- (-2.35,0); 
\draw[thick,gray,densely dotted] (-2.5,.5) -- (-1.4,.9); 
\draw[thick,gray,densely dotted] (-2.35,0) -- (-2.5,-.5); 
\draw[thick,gray,densely dotted] (-2.35,0) -- (-1.35,.25); 
\draw[thick,gray,densely dotted] (-2.35,0) -- (-1.35,-.25); 
\draw[thick,gray,densely dotted] (-2.5,-.5) -- (-2.75,-.8); 
\draw[thick,gray,densely dotted] (-2.5,-.5) -- (-1.4,-.9); 
\draw[thick,gray,densely dotted] (-2.75,-.8) -- (-1.5,-2); 
\draw[thick,dotted] (-2.75,-.8) -- (-3.1,-1.1); 
\draw[thick,dotted] (-1.5,2) -- (-2,2.5); 
\draw[thick,gray,densely dotted] (-1.5,2) -- (-1.4,.9); 
\draw[thick,gray,densely dotted] (-1.5,2) -- (0,2.3); 
\draw[thick,gray,densely dotted] (-1.4,.9) -- (-1.35,.25); 
\draw[thick,gray,densely dotted] (-1.4,.9) -- (0,1.1); 
\draw[thick,gray,densely dotted] (-1.35,.25) -- (-1.1,0); 
\draw[thick,gray,densely dotted] (-1.35,.25) -- (-.825,.4); 
\draw[thick,gray,densely dotted] (-1.1,0) -- (-1.35,-.25); 
\draw[thick,gray,densely dotted] (-1.1,0) -- (0,0); 
\draw[thick,gray,densely dotted] (-1.35,-.25) -- (-1.4,-.9); 
\draw[thick,gray,densely dotted] (-1.35,-.25) -- (-.825,-.4); 
\draw[thick,gray,densely dotted] (-1.4,-.9) -- (-1.5,-2); 
\draw[thick,gray,densely dotted] (-1.4,-.9) -- (0,-1.1); 
\draw[thick,gray,densely dotted] (-1.5,-2) -- (0,-2.3); 
\draw[thick,dotted] (-1.5,-2) -- (-2,-2.5); 
\draw[thick,gray,densely dotted] (0,2.3) -- (0,1.1); 
\draw[thick,gray,densely dotted] (0,2.3) -- (1.5,2); 
\draw[thick,gray,densely dotted] (0,1.1) -- (.6,.7); 
\draw[thick,gray,densely dotted] (0,1.1) -- (-.825,.4); 
\draw[thick,gray,densely dotted] (0,1.1) -- (1.1,1.275); 
\draw[thick,gray,densely dotted] (.6,.7) -- (0,0); 
\draw[thick,gray,densely dotted] (.6,.7) -- (1.9,.75); 
\draw[thick,gray,densely dotted] (-.825,.4) -- (0,0); 
\draw[thick,gray,densely dotted] (0,0) -- (-.825,-.4); 
\draw[thick,gray,densely dotted] (0,0) -- (.6,-.7); 
\draw[thick,gray,densely dotted] (0,0) -- (1.65,0); 
\draw[thick,gray,densely dotted] (-.825,-.4) -- (0,-1.1); 
\draw[thick,gray,densely dotted] (.6,-.7) -- (0,-1.1); 
\draw[thick,gray,densely dotted] (.6,-.7) -- (1.9,-.75); 
\draw[thick,gray,densely dotted] (0,-1.1) -- (0,-2.3); 
\draw[thick,gray,densely dotted] (0,-1.1) -- (1.1,-1.275); 
\draw[thick,gray,densely dotted] (0,-2.3) -- (1.5,-2); 
\draw[thick,dotted] (1.5,2) -- (2,2.5); 
\draw[thick,gray,densely dotted] (1.5,2) -- (2.5,1.1); 
\draw[thick,gray,densely dotted] (1.5,2) -- (1.1,1.275); 
\draw[thick,dotted] (2.5,1.1) -- (3,1); 
\draw[thick,gray,densely dotted] (2.5,1.1) -- (1.9,.75); 
\draw[thick,gray,densely dotted] (1.1,1.275) -- (1.9,.75); 
\draw[thick,gray,densely dotted] (1.9,.75) -- (1.65,0); 
\draw[thick,gray,densely dotted] (1.65,0) -- (1.9,-.75); 
\draw[thick,gray,densely dotted] (1.9,-.75) -- (1.1,-1.275); 
\draw[thick,gray,densely dotted] (1.9,-.75) -- (2.5,-1.1); 
\draw[thick,gray,densely dotted] (1.1,-1.275) -- (1.5,-2); 
\draw[thick,dotted] (2.5,-1.1) -- (3,-1); 
\draw[thick,gray,densely dotted] (2.5,-1.1) -- (1.5,-2); 
\draw[thick,dotted] (1.5,-2) -- (2,-2.5); 

\filldraw[fill=gray,draw=gray] (-2.75,.8) circle (.075); 
\filldraw[fill=gray,draw=gray] (-2.5,.5) circle (.075); 
\filldraw[fill=gray,draw=gray] (-2.35,0) circle (.075); 
\filldraw[fill=gray,draw=gray] (-2.5,-.5) circle (.075); 
\filldraw[fill=gray,draw=gray] (-2.75,-.8) circle (.075); 
\filldraw[fill=gray,draw=gray] (-1.5,2) circle (.075); 
\filldraw[fill=gray,draw=gray] (-1.4,.9) circle (.075); 
\filldraw[fill=gray,draw=gray] (-1.35,.25) circle (.075); 
\filldraw[fill=gray,draw=gray] (-1.1,0) circle (.075); 
\filldraw[fill=gray,draw=gray] (-1.35,-.25) circle (.075); 
\filldraw[fill=gray,draw=gray] (-1.4,-.9) circle (.075); 
\filldraw[fill=gray,draw=gray] (-1.5,-2) circle (.075); 
\filldraw[fill=gray,draw=gray] (0,2.3) circle (.075); 
\filldraw[fill=gray,draw=gray] (0,1.1) circle (.075); 
\filldraw[fill=gray,draw=gray] (.6,.7) circle (.075); 
\filldraw[fill=gray,draw=gray] (-.825,.4) circle (.075); 
\filldraw[fill=gray,draw=gray] (0,0) circle (.075); 
\filldraw[fill=gray,draw=gray] (-.825,-.4) circle (.075); 
\filldraw[fill=gray,draw=gray] (.6,-.7) circle (.075); 
\filldraw[fill=gray,draw=gray] (0,-1.1) circle (.075); 
\filldraw[fill=gray,draw=gray] (0,-2.3) circle (.075); 
\filldraw[fill=gray,draw=gray] (1.5,2) circle (.075); 
\filldraw[fill=gray,draw=gray] (2.5,1.1) circle (.075); 
\filldraw[fill=gray,draw=gray] (1.1,1.275) circle (.075); 
\filldraw[fill=gray,draw=gray] (1.9,.75) circle (.075); 
\filldraw[fill=gray,draw=gray] (1.65,0) circle (.075); 
\filldraw[fill=gray,draw=gray] (1.9,-.75) circle (.075); 
\filldraw[fill=gray,draw=gray] (1.1,-1.275) circle (.075); 
\filldraw[fill=gray,draw=gray] (2.5,-1.1) circle (.075); 
\filldraw[fill=gray,draw=gray] (1.5,-2) circle (.075); 

\end{tikzpicture}
\caption{An arrangement with odd-even invariant zero but no Hamiltonian circuit.  One side of the arrangement, with the tope graph superimposed, is given here.  One tope is lightly shaded for use in a later proof.}
\label{NoHamArr}
\end{figure}

\begin{figure}
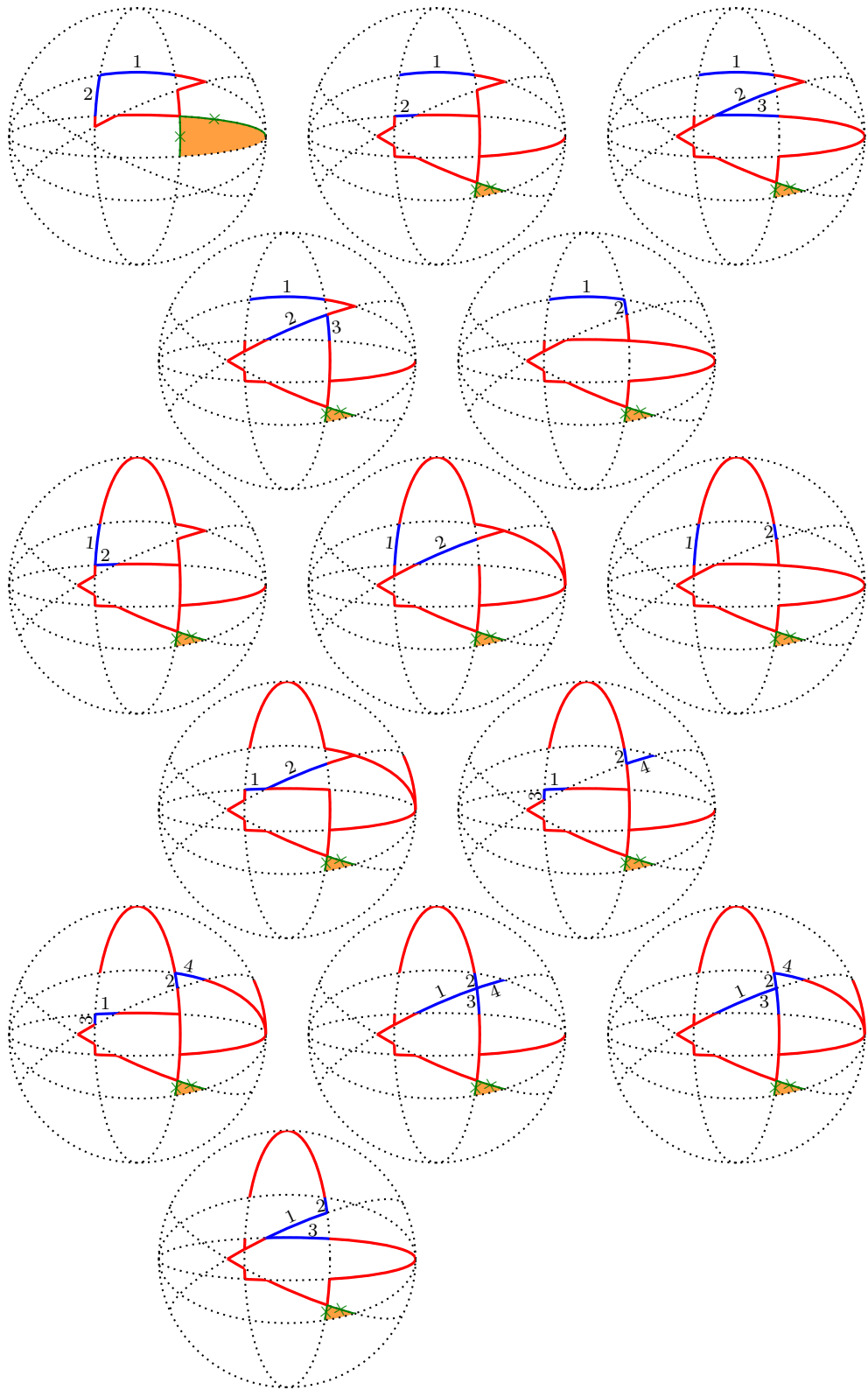

\centering

\caption{An arrangement $\A$ with $\oe(\A)=0$ and non-Hamiltonian tope graph.}
\label{ProofNoHam}
\end{figure}

For a bit of diversity, we give in Figure \ref{NoHamArr} an example of an arrangement $\mathcal{A}$ in $\real^3$ whose tope graph has no Hamiltonian circuit, even though
the odd-even invariant is zero, together with a pictorial proof (Figure \ref{ProofNoHam}).

$\mathcal{A}$ has nine hyperplanes (the eight drawn across the front, and the outer border), thus $\oe(\A)=0$.  We have drawn a part of the tope graph superimposed on the sphere.  If there exists a Hamiltonian circuit, we must be able to select a subgraph of the tope graph such that every vertex has degree two.  This corresponds to selecting a subset $\mathcal{F}$ of the bounding facets of the topes such that for every tope $T$, $|T\cap\mathcal{F}|=2$.  We will show this is impossible, again pictorially, in Figure \ref{ProofNoHam}, by looking at the different possibilities for the upper five-sided, central tope (shaded in Figure \ref{NoHamArr}), and the various associated subcases.  In each case, precisely two facets of the tope must be included.  Numbered, solid (blue) edges indicate the choices we make for a particular case; undecorated, solid (red) edges indicate bounding facets the inclusion of which is forced by those choices; solid, $\times$-ed (green) edges indicate the two bounding facets of a triangular tope the inclusion of which is \textit{not} permitted.  This is the contradiction, as it leaves just one bounding facet that may be selected (and thus, the path through the vertices is unable to continue).  The offending tope in each case is shaded (orange).  Note that there may be more than one ``bad'' tope for a particular set of choices.

\sectb\se.  Bounds on the odd-even invariant.
In the literature there are several
closely related but different settings in which to consider the problems of interest here.
Until now we have considered only the setting of central arrangements in $\real^d$.
By the usual process of \textit{projectivization} (see \cite{ombook}), we may consider instead arrangements in projective space of dimension $d-1$.  When the number of hyperplanes is even, the tope graph for the projective arrangement has odd-even invariant that is half the value of that of the corresponding central arrangement.  When $n$ is odd, the tope graph of the projective arrangement is no longer bipartite and thus the odd-even invariant is not defined
in the analogous way; however, in this case, we may define the odd-even invariant to be zero.

Likewise, by taking the intersections of the hyperplanes of
the central arrangement with a hyperplane missing the origin, we may obtain a non-central
arrangement in $\real^{d-1}$.  This latter process is called \textit{dehomogenization},
the reverse process being \textit{homogenization}, and this setting is of particular importance here, as it allows us to make use of previous research on non-central arrangements, especially non-central planar arrangements.



The question of the maximum ratio $\frac{b}{c}$ for an arrangement of lines in the plane, where $b$ and $c$ are the numbers of topes of each of the two colors, has been studied by several authors; see, e.g., \cite{grunb, palasti, purdywetz, simmwetz}.  It is clear that this maximum ratio must be at least 1, since switching the colors on all the
topes inverts the fraction.  Many have noted that the ratio is less than two.  More precisely, Simmons and Wetzel \cite{simmwetz} showed that $b\le 2c-2-\sum_P (\lambda(P)-2)$, where the sum is over intersections of lines $P$, and $\lambda(P)$ is the number of lines containing $P$.
For centrally simple arrangements in $\real^3$, it is easy to show that $b\le \frac{2}{3}n(n-1)$, with equality if
and only if the triangular regions are precisely those regions that are colored burnt umber.
Thus we see that in the 2-dimensional case, the problem of determining the maximum ratio
is related to that of determining the maximum possible number of
triangles in an arrangement of $n$ lines in the plane. This latter problem has been studied by many authors,
motivated either by a question of Gr\"unbaum in \cite{gruna}, by the 2-coloring problem, or both.  See, for example, \cite{fa, furpal, purdywetz, roud}.  The relatively recent paper of  Bartholdi, Blanc, and Leisel \cite{bbl} gives an account of the problem and adds infinitely many positive integers to the previous infinite list of Forge and Alfons\'in \cite{fa} of integers $n$ for which it is known that there exists an arrangement
of $n$ lines  achieving the bound given by Simmons and Wetzel.

The following theorem shows that there exist arrangements with large odd-even invariant.  This will be used in Section \sf\ to construct further arrangements with odd-even invariant zero and no Hamiltonian circuit.

\thrm\ \te. \begin{enumerate}[(1)]
\item If $d$ is an even positive integer then there is a constant $c_d>0$
such that for arbitrarily large integer values of $n$ there exists an arrangement
$\A$ of $n$ distinct hyperplanes in $\real^d$ for which $\oe(\A) > c_dn^d$.
\item  If $d\ge3$ is odd then there exists $c_d>0$ such that
for arbitrarily large integer values of $n$ there exists an arrangement
$\A$ of distinct hyperplanes in $\real^d$ such that $\oe(\A) > c_dn^{d-1}$.
\item If $d\ge3$ is odd then there exists $\tilde c_d>0$
such that for arbitrarily large integer values of $n$ there exists a {\emph{central}}
arrangement $\A$ of distinct hyperplanes in $\real^d$ for which $\oe(\A) > \tilde c_dn^{d-1}$.
\item If $d\ge4$ is an even positive integer then there exists $\tilde c_d>0$
such that for arbitrarily large integer values of $n$ there exists a \emph{central}
arrangement $\A$ of distinct hyperplanes in $\real^d$ for which $\oe(\A) > \tilde c_dn^{d-2}$.
\end{enumerate}

\proof
The result of Forge and Alfons\'in \cite{fa} previously described implies the validity
for $d=2$, with $c_2 = \frac23$.

If $d = 2m$, $m>1$, we show that $c_d$ may be taken to be $(\frac{2\sqrt{c_2}}{d})^d$.
Assume $m\vert n$, and let $s = \frac n m$.
Let $\A_0 = \{H_1, \ldots, H_s\}$ be an arrangement of $\frac n m$ lines in the plane achieving $\soe(\A_0) > c_2(\frac n m)^2$.
Identify $\real^d$ with $(\real^2)^m$.
For $j = 1, \ldots, m$, let $\pi_j$ be the linear function $\pi_j:\real^d\rightarrow
\real^2$ that takes $(x_1, \ldots, x_{2j-1}, x_{2j}, \ldots,  x_{2m})\in\real^d$
to $(x_{2j-1},x_{2j})\in\real^2$.
For each $(i,j) \in [s]\times [m]$ let $H_{i,j} = \pi_j^{-1}(H_i)$, a hyperplane in $\real^d$,
and let its positive side be the inverse image of the positive side of $H_i$ under $\pi_j$.
Let $\A = \{H_{i,j}: (i,j) \in [s]\times [m]\}$.
The topes of the arrangement $\A$ are the products $T_{k_1}\times \ldots \times T_{k_m}$, where, for each
$r = k_1,\ldots,k_m$, $T_r$ is a tope of $\A_0$.
The signed odd-even invariant of $\A$ is
\begin{eqnarray*}
\soe(\A) & = & \sum_{T_{k_1}, \ldots, T_{k_m}\in \T(\A_0)}(-1)^{\sigma(T_{k_1}\times\ldots\times T_{k_m})}\\
 & = & \sum_{T_{k_1}, \ldots, T_{k_m}\in \T(\A_0)}(-1)^{\sigma(T_{k_1})+\ldots+\sigma(T_{k_m})}\\
 & = & \sum_{T_{k_1}, \ldots, T_{k_m}\in \T(\A_0)}\prod_{j=1}^m(-1)^{\sigma(T_{k_j})}\\
 & = & \prod_{i=1}^m \sum_{T_r\in \T(\A_0)}(-1)^{\sigma(T_r)}\\
 & = & \bigl(\soe(\A_0)\bigr)^m > c_2^m\left(\frac{n}{m}\right)^{2m},
\end{eqnarray*}
where, as before, $\sigma(T)=|\{H:s_T(H)=1\}|$.  From this, (1) follows.

When $d = 2m+1$ with $m\ge1$, we may take $c_d = c_{d-1}$, an appropriate arrangement being
a set of $n$ hyperplanes in $\real^d = \real^{d-1}\times \real$ that are the inverse
images under the projection of $\real^d$ to the first $d-1$ coordinates of
the hyperplanes of an arrangement $\A_0$ in $\real^{d-1}$ that satisfies $\oe(\A_0)
> c_{d-1}n^{d-1}$.

Statement (3) follows from (1) by homogenization, and similarly, statement (4) follows
from (2).
~\endproof

We should note that when the requirement that the hyperplanes be distinct is dropped, statements (1) and (3) are trivial, and in (2) and (4) the order of $n$ can be increased
by one, trivially, by considering arrangements in which each hyperplane appears twice.  In such
an arrangement, the odd-even invariant is the number of topes.

\sectb\sf. Constructions forbidding Hamiltonian circuits.
In this section, we give a theorem that will allow us to build arrangements with odd-even invariant zero and no Hamiltonian circuit.  Before stating this theorem, we give a few necessary definitions.  If $\A$ is an arrangement of $n\geq 1$ hyperplanes and $H\in\A$ then
$\A\setminus\{H\}$ is also an arrangement of $n-1$ hyperplanes in
$\real^d$, called the {\it deletion} of $H$ from $\A$.
Further, we may consider the collection of
intersections $H^\prime\cap H$, where $H^\prime \in \A\setminus\{H\}$,
to be an arrangement of hyperplanes in $\real^{d-1}$.  We
denote this arrangement by $\A/\{H\}$, and call it the {\it restriction} of $\A$ to $H$.
  (The arrangements $\A\setminus\{H\}$
and $\A/\{H\}$ are related to the deletion and contraction operations
of oriented matroids; see \cite{ombook}.)

Given a tope $T$ of $\A$, the restriction of the function $s_T$ to $\A\setminus\{H\}$
 is a tope
of $\A\setminus \{H\}$. If $T$ and $\widetilde T$ are adjacent topes of $\A$ such that the value
of $s_{\widetilde T}$ differs from that of $s_T$ on $H$, then their common restriction
to $\A\setminus\{H\}$ is a tope
of $\A/\{H\}$.

\thrm\ \tf.  Suppose $\A$ is a central arrangement that includes a hyperplane $H$ such that $\oe(\A\setminus\{H\}) > \vert \T(\A/\{H\})\vert$.  Then the tope graph of $\A$ has no Hamiltonian circuit.

\proof
We will say that a tope $T$ of $\A\setminus\{H\}$ {\it straddles} $H$ if $T\cap H^+$ and $T\cap H^-$ are topes of $\A$.  In this case, $T\cap H$ is a tope of $\A/\{H\}$.

Given a $2$-coloring $\gamma_1 : \A \rightarrow \{$chartreuse,
burnt umber$\}$ of the tope graph of $\A$, a $2$-coloring $\gamma_2:\A\setminus\{H\}
\rightarrow \{$chartreuse, burnt umber$\}$ of $\A\setminus\{H\}$ can be obtained
by giving each tope of $\A\setminus\{H\}$ lying in $H^+$ its color as a tope
of $\A$, changing the color of each tope of  $\A\setminus\{H\}$ lying in $H^-$,
and giving each tope $T$ straddling $H$ the color of the tope $T\cap H^+$ of $\A$.

Suppose a Hamiltonian circuit of the tope graph of $\A$ exists.  Remove the
topes (vertices) of the circuit that border $H$ and lie in $H^-$.  The number of
such vertices is $\lvert \T(\A/\{H\})\rvert$, and their removal from the Hamiltonian
circuit leaves a graph consisting of at most $\lvert \T(\A/\{H\})\rvert$ paths.
Each of these paths lies entirely on one side of $H$, so the colors of its vertices
alternate, and the topes along any of these paths contribute at most one to the
odd-even invariant of $\A\setminus\{H\}$.  It follows that $\oe(\A\setminus\{H\})\le
\lvert \T(\A/\{H\})\rvert$.~\endproof

This theorem can be used to verify the nonexistence of Hamiltonian circuits in many
tope graphs.  First, recall from Theorem \te\ that when $d$ is odd there are central arrangements
of $n$ hyperplanes in $\real^d$ having odd-even invariant $\tilde{c}n^{d-1}$, for some positive constant $\tilde{c}$.  On the other hand, it is well-known that in $\real^{d-1}$ any central arrangement has at most $O(n^{d-2})$ topes.  By adding a hyperplane to an arrangement in $\real^d$ that has a large odd-even
invariant, and an even number of hyperplanes, we obtain, according to Theorem \tf, an arrangement with
no Hamiltonian circuit and odd-even invariant zero, as there is an odd number of hyperplanes.  Thus as a corollary of Theorem \te\ and this construction, we have the following theorem.

\thrm\ \th.  If $d\ge3$ is an odd integer there are arbitrarily large central arrangements
for which the odd-even invariant is zero, but whose tope graph admits no
Hamiltonian circuit.

Notice that the number of hyperplanes in an arrangement obtained by the construction
is odd.

For use in the remainder of the paper
 we introduce a refinement of the notation that we have
used.  If $\T$ is any set of topes, we define
\[\soe(\T) = \sum_{T\in\T}(-1)^{\sigma(T)}.\]

If $\A$ is an arrangement of hyperplanes in $\real^d$ then let $\T_b(\A)$ denote the set of bounded topes of $\A$ and let $\T_\infty(\A)$ denote the set of unbounded topes of $\A$.  When $\A$ is a simple arrangement, the topes in $\T_\infty(\A)$ correspond
to the topes of the central arrangement obtained from $\A$ by translating each
hyperplane of $\A$ so that it contains the origin.  We will
indicate this arrangement by $\A_\infty$.  If $\A$ is simple then $\A_\infty$
is a centrally simple arrangement.

\thrm\ \ti.  Suppose $\A$ is a simple arrangement of $n$ hyperplanes in $\real^d$, where
$d$ is odd.  If $n$ is odd, then $\soe(\T(\A))$ $=$ $\soe(\T(A_\infty))$ $=$
 $s\oe(\T(\A_b))$ $= 0$; otherwise, we have $\soe(\T_b(\A))$ $= -\frac{1}{2}\soe(\T(\A_\infty))$.

\proof
Suppose $n$ is odd.
We obtain a centrally simple arrangement of hyperplanes $\A^\prime$ in $\real^{d+1}$ by
viewing the hyperplanes of $\A$ as lying in $\real^d\times \{1\} \subseteq \real^{d+1}$.
We take the hyperplanes of $\A^\prime$ to be those hyperplanes in $\real^{d+1}$
generated by the hyperplanes (in $\real^d\times \{1\}$) of $\A$, together with
 the one additional
hyperplane, $H_0=\real^d\times \{0\}\subseteq\real^{d+1}$,
 whose positive side $H_0^+$ is the open halfspace
containing $\real^d\times \{1\}$.  Since $d+1$ is even and $\A^\prime$ is centrally simple,
 $\soe(\T(\A^\prime)) = 0$, by
Theorem \tc.
The topes of $\A^\prime$ lying in $H_0^+$
are the cones generated by the topes of $\A$,
and the topes of $\A^\prime$ lying in $H_0^-$ are the reflections
of these through the origin.
  Since $n+1$ is even,
 it is clear that $\soe(\T(\A^\prime)) =
-2\soe(\T(\A))$.  Also, clearly, $\soe(\T(\A)) = \soe(\T_b(\A)) + \soe(T_\infty(\A))$.  We have
$\soe(\T_\infty(\A)) = \soe(\T(\A_\infty))$, and, since $n$ is odd, $\soe(\T(\A_\infty)) = 0$.
It follows that $\soe(\T_b(\A)) = 0$ as well, as claimed.

Suppose $n$ is even.
Again we obtain an arrangement $\A^\prime$ in $\real^{d+1}$, this time omitting
the additional hyperplane $H_0$.  Then we have
 $\soe(\T(\A^\prime)) = 2\soe(\T_b(\A)) + \soe(\T(\A_\infty))$.
Since $d+1$ is even and $\A^\prime$ is centrally simple, $\soe(\T(\A^\prime)) = 0$.~\endproof

As a corollary we obtain an asymptotic bound on $\frac b c$, as $n\rightarrow\infty$.

\thrm\ \tj.  Suppose $\A$ is a simple arrangement of $n$ hyperplanes
 in $\real^d$, where $d$ is a fixed odd positive integer, and suppose $\epsilon>0$.
If the tope graph of $\A$ is 2-colored, $b$ topes being colored
burnt umber and $c$, chartreuse, then $\frac b c \le 1+\epsilon $
provided that the number of hyperplanes is sufficiently large.

\proof
We have that $b+c$ is the total number of topes, which is
 $\sum_{k=0}^d \binom n k$ as the arrangement of the $n$ hyperplanes
in $\real^d$ is simple.  This sum is a polynomial of degree $d$ in $n$.
  The odd-even invariant is the difference,
$b-c$, and by Theorem \ti, this is certainly less than the number
of topes of the arrangement $\A_{\infty}$, which is
$2\sum_{k=0}^{\frac{d-1}{2}}\binom{n}{d-1-2k}$, since it is a
centrally simple arrangement.  
This is a polynomial of degree $d-1$ in $n$, and the statement follows.~\endproof

Of course, when $n$ is odd (and assuming $d$ is odd as in the theorem), $\frac b c = 1$ by Theorem \ti.

The following theorem is analogous to Theorem 7, for perfect matchings instead of Hamiltonian
circuits.  It includes the additional hypothesis that the arrangement be centrally simple
and requires a more restrictive bound.

\thrm\ \tk.  Suppose $\A$ is a simple central arrangement of $n$ hyperplanes in $\real^d$ that includes
 a hyperplane $H$ such that $\oe(\A\setminus\{H\}) > 2\lvert \T(\A/\{H\})\rvert$.
Then the tope graph of $\A$ has no perfect matching.

\proof
Since $\oe(\A\setminus\{H\}) \ne 0$, $n-1=\lvert \A\setminus\{H\}\rvert$ is even,
and, by Theorem \tc, $d$ is odd.
Let $b$ be the number of burnt umber topes
and $c$ be the number of chartreuse topes of $\A\setminus\{H\}$ that lie in $H^+$.
We may assume $b\ge c$.
Since $\lvert \A\setminus\{H\}\rvert$ is even, the reflections of topes through
the origin are of the same color.
The remaining topes $T$ of  $\A\setminus\{H\}$ correspond to 
the topes of $\A/\{H\}$, and induce a subgraph of the tope graph of $\A\setminus\{H\}$
that is isomorphic to the tope graph of $\A/\{H\}$.  By Theorem \tc,
$\oe(\A/\{H\}) = 0$.  It follows that the topes 
having a facet in $H$  make no net contribution to the odd-even
invariant of $\A\setminus\{H\}$.
We therefore have $\soe(\A\setminus\{H\}) = 2(b-c)$.

In any perfect matching of the tope graph of $\A$, at least $b-c$ of the topes
lying in $H^+$ and not having $H$ as a facet must be matched to topes in $H^+$
that have $H$ as a facet.  But $b-c$ is larger than $\lvert\T(\A/\{H\})\rvert$,
which is the number of such topes.~\endproof

When $d=3$, Theorem \tk\ and the existence of centrally simple arrangements with
large odd-even invariant establishes the existence of arbitrarily large centrally
simple arrangements in $\real^3$ whose tope graphs have no perfect matchings.

\sectb\sg. Questions.
The foregoing leaves many unanswered questions, of which
we mention only a few.
\begin{enumerate}
\item We have seen that there exist many examples of central arrangements $\A$ with $n$ hyperplanes in $\real^d$ whose tope graphs have no Hamiltonian circuit, even though $\oe(\A) = 0$, with $d$ odd and $n$ odd.
Is the same true for even $d$, or for even $n$?

\item The construction described at the end of Section \sf\ does not yield centrally simple arrangements, in general.  Are there such examples, when the requirement of central simplicity is added?

\item Can the exponent $d-1$ in (4) of Theorem \te\ be increased to $d$?  This would imply a positive answer to Question 1.

\item Is it true that the tope graph of an alternating arrangement $\A$ of $n$ hyperplanes in $\real^d$ has a Hamiltonian circuit if and only if $\oe(\A) = 0$?
\end{enumerate}

\noindent There are analogous questions for non-central arrangements.

\sectc Acknowledgements.
We thank Vic Reiner for his curiosity, which set us off in the first place.  We are grateful to Matthias Beck, and the group at the University of Vienna for many helpful comments and suggestions, and in particular Vivien Ripoll, for alerting us to \cite{conslowil}.  Finally, we thank Javier Bernal for his scrupulous reading of a previous version of the paper and for his many helpful comments, and James Shook for a careful reading of the final version.

\input ref.tex
\vspace{1cm}
\footnotesize
\noindent{\sc Fakult\"at f\"ur Mathematik, Universit\"at Wien, Oskar-Morgenstern-Platz 1, 1090 Wien, Austria}\\
\textit{Email Address:} \texttt{yvonne.kemper@univie.ac.at}\\

\noindent{\sc Department of Mathematical Sciences, George Mason University, Fairfax, VA 22030, U.S.A.}\\
and {\sc National Institute of Standards and Technology, Gaithersburg, MD 20899, U.S.A.}\\
\textit{Email Address:} \texttt{lawrence@gmu.edu}\\

\end{document}

%% file: pre.tex
\usepackage{amssymb,amsmath,amsthm,latexsym}
\long\def\cmnt#1{}

\def\real{{\mathbb R}}

\def\square{{\vcenter{\vbox{\hrule height.4pt
     \hbox {\vrule width.4pt height8pt \kern8pt
       \vrule width.4pt}
     \hrule height.4pt}}}}

\def\proof{\noindent {\sl Proof.\ \ }}
\def\endproof{\hfill$\square$

          \medskip}

\def\abstra#1{{\medskip\narrower\noindent{\bf Abstract.}\sl #1 \medskip}}

\def\sectb#1. #2. {\def\shed{\hbox{\bf #1.}}\penalty -500 \medskip\noindent{\bf #1.\ #2.\ }}
\def\sectc#1. {\penalty -500 \medskip\noindent{\bf #1. }}

\def\prb#1. #2\par{\medskip\noindent{#1. } #2\medskip\par}
\def\prob#1. #2\par{\medskip\noindent{\sc Problem#1. } #2\medskip\par}
\def\thrm#1. #2\par{\medskip\noindent{\sc Theorem#1.} {\sl#2}\medskip\par}
\def\lem#1. #2\par{\medskip\noindent{\sc Lemma#1.} {\sl#2}\par\medskip}
\def\cor#1. #2\par{\medskip\noindent{\sc Corollary#1.} {\sl#2}\medskip\par}
\def\prop#1. #2\par{\medskip\noindent{\sc Proposition#1.} {\sl#2}\medskip\par}
\def\conj#1. #2\par{\medskip\noindent{\sc Conjecture#1.} {\sl#2}\medskip\par}
\def\note#1 {\smallskip\noindent{\bf #1.\enspace}\nobreak}
\def\noteemph#1. #2\par{\medskip\noindent{\sc #1.} {\sl#2}\medskip}
\def\case#1. {\medskip\noindent{\sc #1.} }

\def\oe{\text{\char'33}}
\def\soe{\text{s\char'33}}
\def\T{{\mathcal T}}
\def\A{{\mathcal A}}

\def\sa{1}
\def\sb{2}
\def\sd{3}
\def\se{4}
\def\sf{5}
\def\sg{6}

\def\ta{1}
\def\tb{2}
\def\tc{3}
\def\td{4}
\def\tg{5}
\def\te{6}
\def\tf{7}
\def\th{8}
\def\ti{9}
\def\tj{10}
\def\tk{11}

%% file: hamBrieOrder.bbl
\begin{thebibliography}{99}

\bibitem{bbl}
N. Bartholdi, J. Blanc, and S. Loisel, On simple arrangements of lines
and pseudo-lines in $\mathbb P^2$ and $\real^2$ with the maximum number of triangles.
{\sl Surveys on discrete and computational geometry}, 105--116,
{\bf Contemp. Math.} {\bf 453}, Amer. Math. Soc., Providence, RI, 2008.

\bibitem{ombook}
  A.~Bj\"orner, M.~Las Vergnas, B.~Sturmfels, N.~White, and
G.~M.~Ziegler, {\sl Oriented Matroids}, Encyclopedia of Mathematics,
Cambridge University Press, Cambridge, 1992.

\bibitem{conslowil}
J. H. Conway, N. J. A. Sloane, and A. R. Wilks,
Gray codes for reflection groups.
{\sl Graphs and Combinatorics} {\bf 5} (1989), 315--325.

\bibitem{eag}
R. Eager, A recurrence relation for the odd-even invariant
of bipartite graphs.
{\sl Teknos}, {\bf 10} (Spring, 2001), 33-36.

\bibitem{eaglaw}
R. Eager and J. Lawrence, The odd-even invariant for graphs.
{\sl European Journal of Combinatorics} {\bf 50} (2015), 87--96.

\bibitem{fa}
D. Forge and J. L. Ram\'irez Alfons\'in, Straight line arrangements in the real
projective plane.  {\sl Discrete Comput. Geom.} {\bf 20} (1998), 155--161.

\bibitem{furpal}
Z. F\"uredi and I. Pal\'asti, Arrangements of lines with a large
number of triangles.
{\sl Proc. Amer. Math. Soc.} {\bf 92} (1984), no. 4, 561--566.

\bibitem{gruna}
B. Gr\"unbaum, {\sl Arrangements and Spreads}, Conference Board of the Mathematical
Sciences Regional Conference Series in Mathematics, No. 10, American Mathematical
Society, Providence, R.I., 1972.

\bibitem{grunb}
B. Gr\"unbaum, Two-coloring the faces of arrangements.
{\sl Period. Math. Hungar.} {\bf 11} (1980), no. 3, 181--185.

\bibitem{harb}
H. Harborth, Two colorings of simple arrangements.  In {\sl Finite and
Infinite Sets}, vol. I, II, Egar, 1981, 371--378,
Colloq. Math. Soc. J\'anos Bolyai, 37, 1984.

\bibitem{refoea}
J. Lawrence, The odd-even invariant of an oriented matroid.
{\sl European Journal of Combinatorics} {\bf 27} (2006), 806--813.

\bibitem{palasti}
I. Palasti, The ratio of black and white polygons of a map generated by
general straight lines.
{\sl Period. Math. Hungar.} {\bf 7} (1976), 91--94.

\bibitem{purdywetz}
G. B. Purdy and J. E. Wetzel, Two-coloring inequalities for Euclidean
arrangements in general position.  {\sl Discrete Math.} {\bf 31}
(1980), no. 1, 53--58.

\bibitem{reiner}
V. Reiner, private communication.

\bibitem{roud}
J.-P. Roudneff, On the number of triangles in simple arrangements of
pseudolines in the real projective plane.  {\sl Discrete Mathematics}
{\bf 60} (1986), 243--251.

\bibitem{simmons}
G. J. Simmons, A quadrilateral-free arrangement of sixteen lines.
{\sl Proceedings Amer. Math. Soc.} {\bf 34} (1972), 317--318.

\bibitem{simmwetz}
G. J. Simmons and J. E. Wetzel, A two-coloring inequality for
Euclidean two-arrangements.
{\sl Proc. Amer. Math. Soc.} {\bf 77} (1979), no. 1, 124--127.

\bibitem{stantonwhite}
D. Stanton and D. White, {\sl Constructive Combinatorics},
 Springer, 1986.

\end{thebibliography}
